\documentclass[amstex,12pt,russian,amssymb]{article}

\usepackage{mathtext}
\usepackage[cp1251]{inputenc}
\usepackage[T2A]{fontenc}
\usepackage[russian]{babel}
\usepackage[dvips]{graphicx}
\usepackage{amsmath}
\usepackage{amssymb}
\usepackage{amsxtra}
\usepackage{latexsym}
\usepackage{ifthen}

\renewcommand{\baselinestretch}{1.3}

\begin{document}

\def\Xint#1{\mathchoice
   {\XXint\displaystyle\textstyle{#1}}%
   {\XXint\textstyle\scriptstyle{#1}}%
   {\XXint\scriptstyle\scriptscriptstyle{#1}}%
   {\XXint\scriptscriptstyle\scriptscriptstyle{#1}}%
   \!\int}
\def\XXint#1#2#3{{\setbox0=\hbox{$#1{#2#3}{\int}$}
     \vcenter{\hbox{$#2#3$}}\kern-.5\wd0}}
\def\dashint{\Xint-}

\renewcommand{\abstractname}{}

\noindent
\renewcommand{\baselinestretch}{1.1}
\normalsize \large \noindent

\title{Гельдеровость кольцевых  $Q$-гомеоморфизмов относительно $p$-модуля}

\title{Holder continuity for the ring   $Q$-homeomorphisms with respect  to  $p$-modulus}

\author{ Салимов Р.Р.}

\medskip

 {УДК 517.5}

{\bf  Салимов Р.Р.}  {(Ин-т прикладной математики и механики НАН
Украины, Донецк)}
\medskip

Гельдеровость кольцевых $Q$-гомеоморфизмов относительно $p$-модуля

\begin{abstract} Найдено   достаточное условие гельдеровости  кольцевых $Q$-гомеоморфизмов в $\mathbb{R}^n, n\geq 2$, относительно $p$-модуля
при $n-1<p<n$.

It is founded the sufficient  condition of Holder continuity of  the ring   $Q$-homeomorphisms  in $\mathbb{R}^n, n\geq 2$ with respect  to  $p$-modulus at $n-1<p<n$.

\end{abstract}

{\bf 1. Введение.}  Напомним некоторые определения. Борелева функция
$\varrho:\mathbb{R}^n\to[0,\infty]$ называется {\it допустимой} для
семейства кривых $\Gamma$ в $\mathbb{R}^n$, $n\geq 2$, пишут $\varrho\in{\rm
adm}\,\Gamma$, если
$$\int\limits_{\gamma}\varrho(x)\,ds(x)\geqslant1 $$
для всех $\gamma\in\Gamma$. Пусть $p\geqslant1$. Тогда
$p${\it-модулем} семейства кривых $\Gamma$ называется  величина

$$
M_{p}(\Gamma)=\inf\limits_{\varrho\in{\rm adm}\,\Gamma}
\int\limits_{\mathbb{R}^n}\varrho^{p}(x)\,dm(x).$$
Здесь
$m$ обозначает меру Лебега в $\mathbb{R}^n$.

Пусть $G$ -- область в  $\mathbb{R}^n$, $n\geq 2$. Предположим, что $n-1<p<n$ и \begin{equation}\label{eq1**}
M_p(f\Gamma)\le K\,M_p(\Gamma)\end{equation} для произвольного
семейства $\Gamma$ кривых $\gamma$ в области $G$. При
предположении, что $f$ в (\ref{eq1**}) является гомеоморфизмом,
Герингом  было установлено, что отображение $f$ является {\it липшицевым}, другими словами, при некоторой постоянной $C>0$
и всех $x_0\in G$ справедлива оценка
$$\limsup\limits_{x\to x_0}\frac{|f(x)-f(x_0)|}{|x-x_0|}\leqslant
C,$$ см., напр., теорему 2 в \cite{Ge$_1$}.

\medskip

{\bf 2. О кольцевых $Q$-гомеоморфизмах относительно $p$-модуля.} Всюду далее $B(x_0,\,r)\,=\,\left\{x\in{\Bbb R}^n: \vert
x-x_0\vert\,<\,r\right\},$ ${\Bbb B}^n\,=B(0,1)$,
$\omega_{n-1}$ -- площадь единичной сферы ${\Bbb S}^{n-1}$ в ${\Bbb
R}^n\,.$
Пусть $Q:G\to[0,\infty]$ -- измеримая функция. Для  любого измеримого множества $E\subset \mathbb{R}^n $ и числа $r>0$ обозначим
$$
\dashint\limits_{E}Q(x)\, dm(x)= \frac{1}{m(E)}\int\limits_{E}\,
Q(x)\, dm(x)\,,
$$
$$q_{x_0}(r)=
\frac{1}{\omega_{n-1}r^{n-1}}\int\limits_{S(x_0, r)}Q(x)d\mathcal{A}(x)\,, $$
где $S(x_0, r)=\{x\in
\mathbb{R}^n : |x-x_0|=r\}$, а $d\mathcal{A}(x)$ -- элемент площади
поверхности.

\medskip

Напомним следующие термины. Пусть $E\,,F\,\subseteq\,{\Bbb R}^n\,-$
произвольные множества. Обозначим через $\Delta(E,F;G)$ семейство
всех кривых $\gamma:[a,b]\,\rightarrow\,{\Bbb R^n}\,,$ которые
соединяют $E$ и $F$ в $G\,,$ т.е. $\gamma(a)\,\in\,E\,,\gamma(b)
\,\in\,F$ и $\gamma(t)\,\in\,G\,$ при $a\,<\,t\,<\,b\,.$ Пусть $x_0\in G$,
$d_0=dist (x_0\,,\partial G)$ и
$Q:G\rightarrow\,[0\,,\infty]\,-$ измеримая по Лебегу функция.
Для любых $r_1$ и $r_2$, где $0<r_1<r_2<\infty$, обозначим

$$
R(x_0,r_1,r_2) =\{ x\in {\Bbb R}^n : r_1<|x-x_0|<r_2\},
$$

$$S_{r_i}=S(x_0,r_i)\,,\ \ \
i=1,2.
$$

\medskip Будем говорить, что гомеоморфизм $f:G\to {\Bbb R}^n$ является
{\it кольцевым $Q$-го\-ме\-о\-мор\-физ\-мом относительно $p$-модуля
в точке $x_0\,\in\,G,$} $1<p\leq n$, если соотношение

\begin{equation}\label{2.12}
M_{p}\,\left(\Delta\left(fS_{r_1},fS_{r_2}; fG\right)\right)\leq
\int\limits_{R(x_0,r_1,r_2)} Q(x)\cdot \eta^{p}(|x-x_0|)\ dm(x)
\end{equation}
выполнено для любого кольца $R(x_0,r_1,r_2),$\,\, $0<r_1<r_2<
d_0$ и для каждой измеримой функции $\eta : (r_1,r_2)\to [0,\infty
]\,,$ такой, что
$$
\int\limits_{r_1}^{r_2}\eta(r)\ dr\ \geq\ 1\,.
$$
Говорят, что гомеоморфизм $f:G\to {\Bbb R}^n$ является {\it
кольцевым $Q$-го\-ме\-о\-мор\-физ\-мом относительно $p$-модуля в
области $G$}, если условие $(\ref{2.12})$ выполнено для всех точек
$x_0\,\in\,G\,.$

Развиваемая   в работе теория кольцевых
$Q$-гомеоморфизмов относительно $p$-модуля применима, в частности, к
отображениям квазиконформным в среднем, см. \cite{Go_1}.

\medskip

{\bf 3. Предварительные замечания.}
Следуя работе \cite{MRV}, пару $\mathcal{E}=(A,C)$, где $A\subset\mathbb{R}^n$
-- открытое множество и $C$ -- непустое компактное множество,
содержащееся в $A$, называем {\it конденсатором}. Конденсатор $\mathcal{E}$
называется  {\it кольцевым конденсатором}, если $B=A\setminus C$ --
кольцо, т.е., если $B$ -- область, дополнение которой
$\overline{\mathbb{R}^n}\setminus B$ состоит в точности из двух
компонент.  Говорят
также, что конденсатор $\mathcal{E}=(A,C)$ лежит в области $G$, если $A\subset
G$. Очевидно, что если $f:G\to\mathbb{R}^n$ -- непрерывное, открытое
отображение и $\mathcal{E}=(A,C)$ -- конденсатор в $G$, то $(fA,fC)$ также
конденсатор в $fG$. Далее $f\mathcal{E}=(fA,fC)$.

\medskip
Функция $u:A\to \mathbb{R}$ {\it абсолютно непрерывна на прямой}, имеющей непустое пересечение с $A$, если она абсолютно непрерывна на любом отрезке этой прямой, заключенном в $A$. Функция $u:A\to \mathbb{R}$ принадлежит классу ${\rm ACL}$ ({\it абсолютно непрерывна на почти всех прямых}), если она абсолютно непрерывна на почти всех прямых, параллельных любой координатной оси.

Обозначим через $C_0(A)$  множество
непрерывных функций $u:A\to\mathbb{R}^1$ с компактным носителем,
$W_0(\mathcal{E})=W_0(A,C)$ -- семейство неотрицательных функций
$u:A\to\mathbb{R}^1$ таких, что 1) $u\in C_0(A)$, 2)
$u(x)\geqslant1$ для $x\in C$ и 3) $u$ принадлежит классу ${\rm
ACL}$. Также обозначим
$$
\vert\nabla
u\vert={\left(\sum\limits_{i=1}^n\,{\left(\frac{\partial u}{\partial x_i}\right)}^2
\right)}^{1/2}.$$

При $p\geqslant1$ величину
$${\rm cap_p}\,\mathcal{E}={\rm cap_p}\,(A,C)=\inf\limits_{u\in W_0(\mathcal{E})}\,
\int\limits_{A}\,\vert\nabla u\vert^p\,dm(x)$$ называют
{\it $p$-ёмкостью} конденсатора $\mathcal{E}$. В дальнейшем при  $p>1$ мы
будем использовать равенство
\begin{equation}\label{EMC}
{\rm cap_p}\,\mathcal{E}=M_p(\Delta(\partial A,\partial C; A\setminus C)),\ \
 \end{equation} см. \cite{G}, \cite{H} и
\cite{Sh}.

Известно, что при $1\leq p<n$ \begin{equation}\label{maz} {\rm
cap_p}\,\mathcal{E}\geqslant n{\nu}^{\frac{p}{n}}_n
\left(\frac{n-p}{p-1}\right)^{p-1}\left[m(C)\right]^{\frac{n-p}{n}}\end{equation}
где ${\nu}_n$ -  объем  единичного шара  в ${\Bbb R}^n,\,\,$ см., напр., неравенство (8.9)
в  \cite{Maz}.

При $n-1<p\leq n$ имеет место оценка
 \begin{equation}\label{krd}
\left({\rm cap_p}\,\,
\mathcal{E}\right)^{n-1}\,\ge\,\gamma\,\frac{d(C)^{p}}{m(A)^{1-n+p}}\,\,,
\end{equation}
где $d(C)$ - диаметр компакта $C$,  $\gamma$ - положительная константа, зависящая только от
размерности $n$ и $p\,,$ см. предложение 6  в \cite {Kru}.

\medskip

%SECTION 3 BEGIN
%SECTION 3 BEGIN
%SECTION 3 BEGIN
%SECTION 3 BEGIN
%SECTION 3 BEGIN
%SECTION 3 BEGIN

\medskip

{\bf 4. }{\bf Гёльдеровость кольцевых $Q$-гомеоморфизмов
относительно $p$-модуля.}

Ниже приведена теорема   о  достаточном условии гёльдеровости  кольцевых $Q$-гомеоморфизмов относительно $p$-модуля при
$n-1<p<n$.

\medskip

{\bf Теорема.}  {\it Пусть  $G$ и $G'$ -- области в $\mathbb{R}^n$,
$n\geqslant 2$,    $f:G\to G'$ -- кольцевой $Q$-гомеоморфизм
относительно $p$-модуля {\rm ($n-1<p<n$)}  с $Q(x)\in L_{\alpha}(G)$, $\alpha> \frac{n}{n-p}$ и $F\subset G$ -- произвольный компкт.
Тогда

 \begin{equation}\label{hodler}
|f(x)-f(y)|\leq \lambda_{n,p} \|Q\|_{\alpha}^\frac{1}{n-p} |x-y|^{1-\frac{n}{\alpha(n-p)}}, \ \
\end{equation}
для любой пары точек $x, y\in F$, удовлетворяющих условию $|x-y|<\delta$, где $\delta=\frac{1}{4}dist(F,\partial G)$ и  $\|Q\|_{\alpha}=\left(\int\limits_{G} \, Q^{\alpha}(x)\, dm(x)\right)^{\frac{1}{\alpha}}$ -- норма в пространстве $L_{\alpha}(G) $, $\lambda_{n,p}$ - положительная
постоянная, зависящая только от  $n$ и $p$. }

\medskip

{\it Доказательство.} Рассмотрим сферическое кольцо
$R=R(x,\varepsilon_1, \varepsilon_2)$ с $0<\varepsilon_1<\varepsilon_2<\delta$ такое,
что $R(x,\varepsilon_1, \varepsilon_2)\subset G$. Тогда
$\left(fB\left(x,\varepsilon_2\right),\overline{fB\left(x,\varepsilon_1\right)}\right)$
-- кольцевой конденсатор в   $G'$ и, согласно (\ref{EMC}), имеем
равенство
$$ {\rm cap_p}\ (fB(x,\varepsilon_2),\overline{fB(x
,\varepsilon_1)})=M_{p}(\triangle(\partial
fB(x,\varepsilon_2),\partial f B(x,\varepsilon_1);fR))$$ а ввиду
гомеоморфности  $f,$ равенство
$$\triangle\left(\partial
fB\left(x,\varepsilon_2\right),\partial
fB\left(x,\varepsilon_1\right);fR\right)=f\left(\triangle\left(\partial
B(x,\varepsilon_2) ,\partial
B(x,\varepsilon_1);R\right)\right).$$

Рассмотрим функцию
$$ \eta(t)\,=\,\left
\{\begin{array}{rr} \frac{1}{\varepsilon_2-\varepsilon_1}, &  \ t\in (\varepsilon_1,\varepsilon_2) \\
0, & \ t\in \Bbb{R}\setminus (\varepsilon_1,\varepsilon_2).
\end{array}\right.
$$
В  силу определения  кольцевого  $Q$-гомеоморфизма
  относительно
$p$-модуля, замечаем, что

\begin{equation}\label{eq100} {\rm cap_p}\
(fB(x,\varepsilon_2),\overline{fB(x,\varepsilon_1)}) \leq
\frac{1}{(\varepsilon_2-\varepsilon_1)^p}\int\limits_{R(x,\varepsilon_1,
\varepsilon_2)} Q(x)\ dm(x)\ .\end{equation}

Применяя неравенство Гельдера, имеем
\begin{equation}\label{eq100} {\rm cap_p}\
(fB(x,\varepsilon_2),\overline{fB(x,\varepsilon_1)}) \leq \frac{\left(\Omega_n\varepsilon^n_2\right)^{\frac{\alpha-1}{\alpha}}}{(\varepsilon_2-\varepsilon_1)^p}
\|Q\|_{\alpha} .\end{equation}

Далее, выбирая $\varepsilon=|x-y|$,
$\varepsilon_1=2\varepsilon$ и $\varepsilon_2=4\varepsilon$, получим

\begin{equation}\label{eq101}{\rm cap_p}\ (fB(x,4\varepsilon),f\overline{B(x,2\varepsilon)})\le\,\gamma_1\|Q\|_{\alpha}\cdot\varepsilon^{\frac{\alpha n-\alpha p-n}{\alpha}}\,.
\end{equation}
С другой стороны,  в силу  неравенства (\ref{maz}) вытекает оценка

\begin{equation}\label{eq102} {\rm cap_p}\ (fB(x,4\varepsilon),f\overline{B(x,2\varepsilon)})
\ge \gamma_2\left[m(fB(x,2\varepsilon))\right]^{\frac{n-p}{n}}\,,
\end{equation}
где    $\gamma_2$  -- положительная константа, зависящая только от размерности
пространства   $n$ и $p.$

Комбинируя   (\ref{eq101}) и (\ref{eq102}), получаем, что
\begin{equation}\label{eq4.2} m(fB(x,2\varepsilon)) \leqslant \gamma_3 \|Q\|^{\frac{n}{n-p}}_{\alpha}\varepsilon^{\frac{n(\alpha n-\alpha p-n)}{\alpha(n-p)}} \,,\end{equation} где $\gamma_3$ -
положительная постоянная зависящая только от $n$ и $p$.

Далее, выбирая в (\ref{eq100}) $\varepsilon_1=\varepsilon$ и
$\varepsilon_2=2\varepsilon$, получим
\begin{equation}\label{eq91}{\rm cap_p}\ (fB(x,2\varepsilon),f\overline{B(x,\varepsilon)})\le\,
\gamma_4\|Q\|_{\alpha}\cdot\varepsilon^{\frac{\alpha n-\alpha p-n}{\alpha}}\,.
\end{equation}
С другой стороны, в силу неравенства (\ref{krd}), получаем
\begin{equation}\label{eq10*!} {\rm cap_p}\ (fB(x,2\varepsilon),f\overline{B(x,\varepsilon)})
\ge
\left(\gamma_5\frac{d^p(\overline{fB(x,\varepsilon)})}{m^{1-n+p}(fB(x,2\varepsilon))}\right)^{\frac{1}{n-1}}\,,
\end{equation}
где    $\gamma_5$  --   положительная константа, зависящая только от  $n$ и $p.$

Комбинируя   (\ref{eq91}) и (\ref{eq10*!}), получаем, что
$$
d(\overline{fB(x,\varepsilon)}) \le \gamma \|Q\|_{\alpha}^\frac{1}{n-p} \varepsilon^{1-\frac{n}{\alpha(n-p)}} \,.$$
где    $\gamma$  --   положительная константа, зависящая только от  $n$ и $p.$ Оценка (\ref{hodler})  получается отсюда и из очевидного неравенства $d(\overline{fB(x,\varepsilon)})\geq |f(x)-f(y)|\,.$

\medskip

Салимов Руслан Радикович

Институт прикладной математики и механики НАН Украины

ул. Розы  Люксембург 74, Донецк, 83114.

Рабочий телефон: 311-01-45

Email: salimov07@rambler.ru, ruslan623@yandex.ru,

\end{document}